\documentclass[11pt,a4paper]{article}
\usepackage{amsfonts,here,graphicx}
\usepackage{cite,amsfonts,amssymb} 

\def\proof{\noindent \medskip {\bf Proof}$\;\;$}

\newtheorem{prop}{Proposition}

\newtheorem{lem}{Lemma}

\newtheorem{theor}{Theorem}

\newtheorem{defin}{Definition}

\begin{document}

\begin{center}
\Large{\bf STOCHASTIC VOLTERRA EQUATIONS OF NONSCALAR TYPE
IN HILBERT SPACE\\[25pt]}
\end{center}

\noindent{\bf Karczewska, Anna\\[5pt]}
{\it Department of  Mathematics,
University of Zielona G\'ora,  
ul.\ Szafrana 4a, 65-246 Zielona G\'ora, Poland
{\bf [A.Karczewska@im.uz.zgora.pl]}\\}

\begin{center}
{\bf Abstract}
\end{center}

{\small
\noindent
In the paper stochastic Volterra equations of nonscalar type
in Hilbert space are studied. 
The aim of the paper is to provide some results on stochastic
convolution and mild  solutions to those Volterra equations.
The motivation of the paper comes 
from a model of aging viscoelastic materials.
The pseudo-resolvent approach is used.\\[10pt]
}


\section{Introduction}\label{intro}

Assume that $H$ is a separable Hilbert space with a scalar product
$\langle\cdot ,\cdot\rangle_H$ and a norm $|\cdot |_H$.
We are concerned with a stochastic linear Volterra equation
of nonscalar type in $H$ of the form
\begin{equation} \label{deq1}
 X(t) = X_0 + \int_0^t A(t-\tau) X(\tau) d\tau 
 +\int_0^t  \Psi(\tau)dW(\tau), \quad t\in [0,T].
\end{equation}

Let $G$ denote another separable Hilbert space such that 
$G\stackrel{d}{\hookrightarrow} H$ ($G$ is densely imbedded in $H$).
In (\ref{deq1}), 
$A\in L^1_\mathrm{loc}(\mathbb{R}_+;\mathcal{B}(G,H))$,
where $\mathcal{B}(G,H)$ denotes the space of all bounded linear 
operators acting from $G$ into $H$. 
We assume that $W$ is a cylindrical Wiener process with values in some
another Hilbert space $U$, and $\Psi$ is a stochastic process specified
below.

Observe that Volterra equation of so-called scalar type 
(see (Pr\"uss 1993))
corresponds to the operator $A(t)$ of the form
 $A(t)=a(t)A$, where 
$a\in L^1_\mathrm{loc}(\mathbb{R}_+)$, and $A$ is a closed linear 
unbounded operator densely defined in $H$ and $G=(D(A),|\cdot|_A)$,
where $D(A)$ denotes the domain of $A$ equipped with the graph norm
$|\cdot|_A$ of $A$.

The setting for (\ref{deq1}) when $A(t)$ are bounded 
seems to cover many abstract treatments of Volterra equations
although this is not the most general approach.
However, it is very difficult to develop a reasonable theory 
in very general case even for deterministic equation.  

As we have already written, the paper has been inspired by a model
of a hereditarily-elastic anisotropic aging body with a straight crack,
see (Costabel et al.\ 2004).

The equation (\ref{deq1}) is a stochastic version of the deterministic
Volterra equation 
\begin{equation} \label{deq2}
 u(t) =  \int_0^t A(t-\tau) u(\tau) d\tau  + f(t)\,, \quad t\geq 0\,,
\end{equation}
studied extensively in (Pr\"uss 1993).

\begin{defin}
A family $(\widetilde{S}(t))_{t\geq 0}\subset \mathcal{B}(H)$ is called 
{\tt pseudo-resolvent} for (\ref{deq2}) if the following conditions are
satisfied:
\begin{description}
\item[(S1)] $\widetilde{S}(t)$ is strongly continuous in $H$ on 
 $\mathbb{R}_+$, and $\widetilde{S}(0)=I$;
\item[(S2)] $U(t)=\int_0^t \widetilde{S}(\tau)d\tau$ is leaving $G$ invariant,
 and $(U(t))_{t\geq 0}\subset \mathcal{B}(G)$ is locally Lipschitz on
  $\mathbb{R}_+$;
\item[(S3)] the following resolvent equations hold
\begin{eqnarray} 
 \widetilde{S}(t)y &=& y+ \int_0^t A(t-\tau) dU(\tau)y,
 \quad t\geq 0,\quad y\in G\;,
 \label{deq3} \\
 \widetilde{S}(t)y &=& y+ \int_0^t \widetilde{S}(t-\tau) 
 A(\tau)y\,d\tau,\quad t\geq 0,\quad y\in G\;.
 \label{deq3a} 
\end{eqnarray}
\begin{itemize} 
\item Equations (\ref{deq3}) and (\ref{deq3a}) are called the\/ 
{\tt first}, resp. {\tt second resol\-vent equation} for  (\ref{deq2}).
\item A pseudo-resolvent $\widetilde{S}(t)$ is called {\tt resolvent} 
for (\ref{deq2})
if in addidion 
\end{itemize}
\item[(S4)] for $~y\in G,~~\widetilde{S}(\cdot)y\in G$ a.e.\ and 
$\widetilde{S}(\cdot)y$ is
Bochner-measurable in $G$ on $\mathbb{R}_+$.
\end{description}
\end{defin}
{\bf Comment }{ \em
Because $A(t)\in\mathcal{B}(G,H)$ the convolution $A\star dU$ in
(\ref{deq3}) is well-defined though the function $A(\cdot)y$ is not
assumed to be continuous. Indeed, for $y\in G$ the function
$g(t) := U(t)y$ is locally Lipschitz in $G$ by the condition 
{\bf (S2)} and then $g\in BV_\mathrm{loc}(\mathbb{R}_+;G)$.
For details, see (Pr\"uss 1993, Section~6).}

Let us emphasize that the pseudo-resolvents are always unique. We shall 
assume in the paper that the equation (\ref{deq1}) is 
{\tt well-posed} in this sense that (\ref{deq1}) admits a 
pseudo-resolvent $\widetilde{S}(t),~t\geq 0$. 
Precise definition of well-posedness
is given in (Pr\"uss 1993). That definition is a direct extension of
well-posedness of Cauchy problems and
Volterra equations of scalar type. The lack of well-posedness
leads to distribution resolvents, see e.g.\ (DaPrato and Iannelli 1984).

In sections \ref{prob} and \ref{prop} we shall use the following
{\sc Volterra Assumptions} (abbr.~(VA)): 
{\it \begin{enumerate} 
\item $A\in L_\mathrm{loc}^1 (\mathbb{R}_+;\mathcal{B}(G,H))$;
\item  $\widetilde{S}(t),~t\geq 0$, are pseudo-resolvent operators 
 corresponding to (\ref{deq2}). 
\end{enumerate}
}

\section{Probabilistic background}\label{prob} 

We are given a probability space 
$(\Omega,\mathcal{F},(\mathcal{F}_t),P),~t\geq 0$, 
with normal filtration and a Wiener process $W$ with the positive 
linear covariance operator $Q$ in $U$. Assume that  
$W$ is a cylindrical process, that is 
 $\mathrm{Tr}\,Q=+\infty$ (if $\mathrm{Tr}\,Q<+\infty$ 
then $W$ is a genuine Wiener process). 
We will need the subspace $U_0:=Q^{1/2}(U)$ of the space $U$, which 
endowed with the inner product 
$\langle u,v\rangle_{U_0}:=\langle Q^{-1/2}u, Q^{-1/2}v\rangle_U$
forms a Hilbert space. 

A linear, bounded operator $B$ acting from $U_0$
into $H$ is called a {\em Hilbert-Schmidt} if
$\sum_{k=1}^{+\infty} |Bu_k|_H^2 < +\infty $, where $\{u_k\}\subset U_0$
is a base in $U_0$. The set $L_2(U_0,H)$ of all Hilbert-Schmidt operators
from $U_0$ into $H$, equipped with the norm
$|B|_{L_2(U_0,H)}:=(\sum_{k=1}^{+\infty} |Bu_k|_H^2)^{1/2}$, 
is  a separable Hilbert space.
For abbreviation we shall denote  $L_2^0:=L_2(U_0,H)$.

Assume that $\Psi$ belongs to the 
class of $L_2^0$-predictable processes satisfying condition
$$
 P\left(\int_0^T | \Psi(\tau) |_{L_2^0}^2 \,d\tau < +\infty\right) = 1.
$$
Such processes are called {\em stochastically integrable} on $[0,T]$
and create a linear space which will be denoted by 
$\mathcal{N}(0,T;L_2^0)$. 

Let $\Phi(t),~t\in[0,T]$, be a measurable $L_2^0$-valued process. 
Define the norms
\begin{eqnarray*}
 ||\Phi||_t &:=& \left\{\mathbb{E}\left( \int_0^t |\Phi(\tau)|_{L_2^0}^2\,d\tau 
 \right) \right\}^{\frac{1}{2}}  
 = \left\{\mathbb{E} \int_0^t 
 \left[ \mathrm{Tr} (\Phi(\tau)Q^{\frac{1}{2}}) (\Phi(\tau)Q^{\frac{1}{2}})^*
 \right] d\tau \right\}^{\frac{1}{2}}, ~~t\in [0,T].
\end{eqnarray*}

 By $\mathcal{N}^2(0,T;L_2^0)$ we shall denote a Hilbert space of 
all $L_2^0$-predictable processes $\Phi$ such that $||\Phi ||_T <+\infty$. \\

In the whole paper we shall use the following {\sc Probability
Assumptions} (abbr.~(PA)):
{\it \begin{enumerate} 
\item $X_0$ is an $H$-valued, $\mathcal{F}_0$-measurable random variable;
\item  $\Psi\in \mathcal{N}^2(0,T;L_2^0)$, where the interval $[0,T]$ is fixed. 
\end{enumerate}
}

\begin{defin} \label{def5}
Assume that conditions (VA) and (PA) hold.
An $H$-valued predictable process $X(t),~t\in [0,T]$, is said to be a 
{\tt strong solution} to  (\ref{deq1}), if $X$ has a version such that
$P(X(t)\in G)=1$ for almost all $t\in [0,T]$,~
$\int_0^T |A(t-\tau)X(\tau)|_H d\tau<+\infty,~ P-a.s.$ and for any 
$t\in [0,T]$ the equation (\ref{deq1}) holds $P$-a.s.
\end{defin}
\begin{defin} \label{def6}
Let conditions (VA) and (PA) hold.
An $H$-valued predictable process $X(t),~t\in [0,T]$, is said to be a
{\tt mild solution} to the stochastic Volterra equation (\ref{deq1}), if
\begin{equation} \label{deq4_0}
 \mathbb{E}\left( \int_0^t |
 \widetilde{S}(t-\tau)\Psi(\tau)|_{L_2^0}^2 \,d\tau\right)
 \leq + \infty \quad \hbox{for} \quad t\leq T
\end{equation}
and, for arbitrary $t\in [0,T]$,
\begin{equation}\label{deq4}
 X(t) = \widetilde{S}(t)X_0 + \int_0^t
\widetilde{S}(t-\tau)\Psi(\tau)\,dW(\tau),
\end{equation}
where $\widetilde{S}(t)$ is the pseudo-resolvent for the equation (\ref{deq2}). 
\end{defin}

\noindent{\bf Comment }{\em The well-posedness of (\ref{deq1})
implies existence and uniqueness of the resolvent $\widetilde{S}(t),~t\geq 0$,
and in the consequence, existence and uniqueness of the 
mild solution to the equation (\ref{deq1}).}

\begin{prop} \label{dl1} 
 Assume that (VA) hold, $B$  
 is a linear bounded operator acting from $U$ into the space $H$ and
\begin{equation} \label{deq7} 
  \int_0^T |\widetilde{S}(\tau)B|_{L_2^0}^2\,d\tau =  \int_0^T \mathrm{Tr} 
  [\widetilde{S}(\tau) BQB^*\widetilde{S}^*(\tau)]\,d\tau < +\infty,
\end{equation} 
where $B^*, \widetilde{S}^*(\tau)$ are appropriate adjoint operators.

Then
\begin{description}
 \item[~~(i)] the process $\widetilde{W}^B:=\int_0^t 
 \widetilde{S}(t-\tau)BdW(\tau)$ 
 is Gaussian, mean-square continuous on [0,T] 
 and then has a predictable version;
 \item[~(ii)] 
\begin{equation} \label{deq8} 
  \mathrm{Cov}~\widetilde{W}^B(t) =  \int_0^t [\widetilde{S}(\tau) 
  BQB^*\widetilde{S}^*(\tau)]\,d\tau,   \quad   t\in [0,T];
\end{equation} 
 \item[(iii)] trajectories of the process $\widetilde{W}^B$ are P-a.s.\ 
 square integrable on [0,T].
\end{description}
\end{prop}

\noindent{\bf Comment }{ \em
Proposition \ref{dl1} is analogous to the known result obtained in
the semigroup case. Gaussianity of the process $\widetilde{W}^B$ follows from
the definition and properties of stochastic integral. 
Part {\bf (ii)} comes from theory of stochastic integral while 
{\bf (iii)} follows from the definition of $\widetilde{W}^B$ and the assumption
(\ref{deq7}).  }

\begin{prop} \label{dl2}
Assume that conditions (VA) hold with operator 
$A\in W^{1,1}([0,T];\mathcal{B}(G,H))$. 
Let $X$ be a strong solution to the equation (\ref{deq1}) in the
case $\Psi(t)=B$, where $B\in \mathcal{B}(U,H)$ and trajectories of $X$ are
integrable
w.p.~1 on $[0,T]$.  
Then, for any function $\xi\in C^1([0,T];H^*)$ and
$t\in [0,T]$, the following formula holds
\begin{eqnarray} \label{deq10} 
  \langle X(t),\xi(t)\rangle_H & = &  \langle X_0,\xi(0)\rangle_H
 + \int_0^t  \langle  (\dot{A}\star X)(\tau)
  + A(0)X(\tau),\xi(\tau)\rangle_H  d\tau  \nonumber \\ & + &
  \int_0^t \langle \xi(\tau),BdW(\tau)\rangle_H  
  +\int_0^t \langle X(\tau),\dot{\xi}(\tau)\rangle_H d\tau,
\end{eqnarray} 
where dots above $A$ and $\xi$ mean time derivatives and $\;\star\,$ 
means the convolution.
\end{prop}

\proof{First, we consider functions  
$\xi(\tau):=\xi_0\varphi(\tau)$,
$\tau\in [0,T]$, where $\xi_0\in H^*$ and $\varphi\in C^1[0,T]$.
Denote 
$ F_{\xi_0}(t) := \langle X(t),\xi_0 \rangle_H , \quad t\in [0,T].$
\\[1mm]
Using It\^o's formula to the process $ F_{\xi_0}(t)\varphi(t)$, we obtain
\begin{equation} \label{deq11} 
 d[F_{\xi_0}(t)\varphi(t)] = \varphi(t)dF_{\xi_0}(t) 
 + \dot{\varphi}(t) F_{\xi_0}(t)dt, \quad\quad t\in[0,T].
\end{equation} 
Then 
\begin{eqnarray} \label{deq12} 
 dF_{\xi_0}(t) & = &\langle\int_0^t \dot{A}(t-\tau)X(\tau)d\tau + A(0)X(t),
 \xi_0\rangle_H dt +\langle BdW(t),\xi_0\rangle_H \nonumber \\
 & = & \langle (\dot{A}\star X)(t) +A(0)X(t),\xi_0\rangle_H dt + 
 \langle BdW(t),\xi_0 \rangle_H .
\end{eqnarray} 
From (\ref{deq11}) and (\ref{deq12}), 
\begin{eqnarray*}  
 F_{\xi_0}(t)\varphi(t) & = & F_{\xi_0}(0)\varphi(0) +
 \int_0^t \varphi(s) \langle (\dot{A}\star X)(s) 
 + A(0)X(s), \xi_0\rangle_H ds \\ && + \int_0^t 
 \langle\varphi(s)BdW(s),\xi_0\rangle_H + \int_0^t \dot{\varphi}(s)
 \langle X(s), \xi_0\rangle_H ds \\
  & = & \langle X_0,\xi(0)\rangle + \int_0^t \langle (\dot{A}\star X)(s)
 + A(0)X(s), \xi(s)\rangle_H ds \\ 
 && +\int_0^t \langle BdW(s),\xi(s)\rangle_H 
 +\int_0^t \langle X(s),\dot{\xi}(s)\rangle_H ds .
\end{eqnarray*} 
Hence,~ we proved the formula ~(\ref{deq10})~ for functions~ $\xi$~ of the
form~ 
$\xi(s)=\xi_0\varphi(s)$,~ $s\in [0,T]$.~
 Because such functions form a dense 
subspace in~ $C^1([0,T];H^*)$,~ the proposition is true.
~~\hspace{2ex}\hfill  ~$\blacksquare$}

\section{Properties of stochastic convolution}\label{prop}

In this section we study mild solution to the equation 
(\ref{deq1}). We use the following notation 
\begin{equation}\label{deq5}
 \widetilde{W}^\Psi(t) :=  \int_0^t \widetilde{S}(t-\tau)\Psi(\tau)\,
 dW(\tau),\quad t\in [0,T],
\end{equation}
where the condition (\ref{deq4_0}) holds.

\begin{prop} \label{pr3}
Assume that $\widetilde{S}(t),~t\geq 0$, are the pseudo-resolvent operators 
to the Volterra equation (\ref{deq2}). Then, for arbitrary  
$\Psi\in\mathcal{N}(0,T;L_2^0)$, the process $\widetilde{W}^\Psi(t),~t\geq 0$, 
given by (\ref{deq5}) has a predictable version.
\end{prop}

\proof{Proof is analogous to construction  
of stochastic integral, see e.g.\ (Liptser and Shiryayev 1973).

The process $\widetilde{S}(t-\tau)\Psi(\tau),~\tau\in[0,T]$,
belongs to $\mathcal{N}(0,T;L_2^0)$, because $\Psi\in\mathcal{N}(0,T;L_2^0)$.
We may use the following estimate:  
for arbitrary $a>0,~b>0$ and $t\in[0,T]$, 
\begin{equation} \label{g1}
P(|\widetilde{W}^\Psi(t)|_H>a) \leq \frac{b}{a^2}+P\left( \int_0^t  
 |\widetilde{S}(t-\tau)\Psi(\tau)|_{L_2^0}^2 d\tau >b\right)\;.
\end{equation}
Since pseudo-resolvent operators $\widetilde{S}(t),~t\geq 0$, are uniformly
bounded on compact itervals, there exists a constant $C>0$, such that 
$|\widetilde{S}(t-\tau)\Psi(\tau)|_{L_2^0}^2\leq C^2|\Psi(\tau)|_{L_2^0}^2$,
$\tau\in[0,T]$.

Then (\ref{g1}) reads
\begin{equation} \label{g2}
P(|\widetilde{W}^\Psi(t)|_H>a) \leq \frac{b}{a^2}+P\left( \int_0^t  
|\Psi(\tau)|_{L_2^0}^2 d\tau >\frac{b}{C^2}\right)\;.
\end{equation}

We prove predictability of the process 
$\widetilde{W}^\Psi$ in two steps. 
In the first step $\Psi$ is an elementary process, so the process
$\widetilde{W}^\Psi$ has a predictable 
version by Proposition \ref{dl1}, part~{\bf (i)}.

In the second step $\Psi\in\mathcal{N}(0,T;L_2^0)$. There exists a sequence 
$(\Psi_n)$ of elementary processes that for arbitrary $c>0$,
\begin{equation} \label{g3}
P\left( \int_0^T |\Psi(\tau)-\Psi_n(\tau)|_{L_2^0}^2 d\tau >c \right)~
\longrightarrow_{\hspace{-4.5ex}_{n\to +\infty}} 0\;.
\end{equation}
Because, by the previous part of the proof, the sequence $(W_n^\Psi)$
converges in probability, it has subsequence converging almost surely.
This fact implies the predictability of $W^\Psi(t),~t\in[0,T]$.
~\hfill $\blacksquare$}
\begin{prop} \label{pr3a}
Assume that $\Psi\in\mathcal{N}^2(0,T;L_2^0)$. Then the process 
$~\widetilde{W}^\Psi(t)$, $t\geq 0$, defined by (\ref{deq5})
 has square integrable 
trajectories.
\end{prop}
\proof{
From Fubini's theorem and boundness of operators $\widetilde{S}(t)$ we obtain
\begin{eqnarray*}
 \mathbb{E}\int_0^T \!\left|\! \int_0^t\!
\widetilde{S}(t-\tau)\Psi(\tau)dW(\tau)\right|_H^2 dt
  \!&\! =\!&\! \int_0^T \! \left[\mathbb{E}\left| 
 \int_0^t \widetilde{S}(t-\tau)\Psi(\tau)dW(\tau)\right|_H^2\right] dt \\
 = \int_0^T\int_0^t |\widetilde{S}(t-\tau)\Psi(\tau)|_{L_2^0}^2\; d\tau dt 
 \!&\! \leq \!&\!
 M \int_0^T\int_0^t |\Psi(\tau)|_{L_2^0}^2 \; d\tau dt  
 < +\infty, \quad  t\in [0,T].  
\end{eqnarray*}
\hfill $\blacksquare$}

\begin{theor} \label{pr4}
Assume that (VA) and (PA) hold with $A\in W^{1,1}([0,T],\mathcal{B}(G,H))$, 
and $X$ is an $H$-valued predictable process. Let $S(t)$ be a weak
analytic resolvent for (\ref{deq1}). 
Then strong solution to (\ref{deq1}) is a mild solution, that is
\begin{equation}\label{deq9a}
X(t)= X_0+\widetilde{W}^\Psi(t) \quad t\in [0,T].
\end{equation}
\end{theor}

\proof{
By Proposition \ref{dl2}, for any 
$\xi\in C^1([0,T],H^*)$ and $t\in [0,T]$, the following equation holds
\begin{eqnarray}\label{deq14}
 \langle X(t),\xi(t)\rangle_H &\!=\!& \langle X_0,\xi(0)\rangle_H  +
\int_0^t\langle (\dot{A}\star X)(\tau) +
 A(0)X(\tau),\xi(\tau)\rangle_H\, d\tau  \\
 &\!+\!& \int_0^t \langle \Psi(\tau) dW(\tau),\xi(\tau)\rangle_H  
 + \int_0^t \langle X(\tau),\dot{\xi}(\tau)  \rangle_H\, d\tau, 
 ~\mathrm{~P-a.s.} \nonumber
\end{eqnarray}
Let us take $\xi(\tau):=\widetilde{S}^*(t-\tau)\zeta$, 
for $\zeta\in H^*$, 
 $\tau\in [0,t]$.

 Now, the equation (\ref{deq14}) may be written like
\begin{eqnarray*}  
 \langle X(t),\widetilde{S}^*(0)\zeta\rangle_H & = & 
 \langle X_0,\widetilde{S}^*(t)\zeta\rangle_H +
 \int_0^t\langle (\dot{A}\star X)(\tau)
 + A(0)X(\tau),\widetilde{S}^*(t-\tau)\zeta\rangle_H\, d\tau  \\
 &+& \int_0^t \langle \Psi(\tau) dW(\tau),\widetilde{S}^*(t-\tau)
  \zeta\rangle_H 
 + \int_0^t \langle X(\tau),(\widetilde{S}^*(t-\tau)\zeta)'\rangle_H\, 
 d\tau, 
\end{eqnarray*}
where derivative ()' in the last term is taken over $\tau$. 

Since $\widetilde{S}^*(0)=I$, we obtain
\begin{eqnarray}\label{deq14a}
 \hspace{-4ex} \langle X(t),\zeta\rangle_H & = &
 \langle S(t)X_0,\zeta\rangle_H +
 \int_0^t \!\langle \widetilde{S}(t-\tau) \left[\int_0^\tau \!\!
 \dot{A}(\tau-\sigma)X(\sigma)d\sigma 
 +A(0)X(\tau)\right],\zeta\rangle_H\, d\tau~~ \nonumber\\
 & + &\int_0^t \langle \widetilde{S}(t-\tau)\Psi(\tau)dW(\tau),
 \zeta\rangle_H +
 \int_0^t \langle \dot{\widetilde{S}}(t-\tau)X(\tau),
 \zeta\rangle_H\, d\tau\, 
\end{eqnarray}
for any $\zeta\in H^*$.

Let us analyze the right hand side of (\ref{deq14a}).
From the properties of convolution and {\bf (S3)},
\begin{eqnarray*}  
 \int_0^t \dot{\widetilde{S}}(t-\tau)X(\tau) d\tau \!&\!=\!&\! 
 -\int_0^t \dot{\widetilde{S}}(\tau)X(t-\tau) d\tau 
  = -\int_0^t\!\left[ \!\int_0^\tau\! 
A(\tau-s)dU(s)\right]'\!X(t-\tau)d\tau \\
 \!&\!=\!&\! -\int_0^t\left[ 
 \int_0^\tau \dot{A}(\tau-s)dU(s)+A(0)\widetilde{S}(\tau)\right]
 X(t-\tau) d\tau.
\end{eqnarray*}
We have 
$$ \int_0^t  A(0)\widetilde{S}(t-\tau)X(\tau)\, d\tau =
  \int_0^t A(0)\widetilde{S}(\tau)X(t-\tau),\, d\tau $$
 and, from  {\bf (S2)},
$$
 \int_0^\tau \dot{A}(\tau-s)dU(s) =\int_0^\tau\dot{A}(\tau-s)U'(s)ds =
 \int_0^\tau \dot{A}(\tau-s)\widetilde{S}(s)ds.
$$ 
Hence
\begin{eqnarray*} 
&& \int_0^t \widetilde{S}(t-\tau) \left[\int_0^\tau
 \dot{A}(\tau-\sigma)X(\sigma)d\sigma\right] d\tau =
 \int_0^t \widetilde{S}(t-\tau)(\dot{A}\star X)(\tau) d\tau = \\[2mm]
 &\!\!=\!\!& (\widetilde{S}\star (\dot{A}\star X)(\tau))(t) =
 ((\widetilde{S}\star\dot{A})(\tau)\star X)(t) = \\[2mm]
 &\!\!=\!\!&  \int_0^t (\widetilde{S}\star\dot{A})(\tau)X(t-\tau)d\tau
 = \int_0^t \left[ \int_0^\tau \widetilde{S}(\tau)\dot{A}(\tau-s)ds\right]
 X(t-\tau) d\tau\,. 
\end{eqnarray*} 
So, the right hand side of (\ref{deq14a}) reduces to
$\int_0^t \langle \widetilde{S}(t-\tau)\Psi(\tau)dW(\tau),\zeta
\rangle_H$, what means that (\ref{deq9a}) holds.
\hfill $\blacksquare$}\\

\begin{theor} \label{th1}
Assume that $A\in L^1_\mathrm{loc}([0,T];\mathcal{B}(G,H))$ and
$\widetilde{S}(t),~t\geq 0$, is a resolvent to the equation (\ref{deq2}).
If $\Psi\in \mathcal{N}^2(0,T;L_2^0)$, then the stochastic
convolution $\widetilde{W}^\Psi$ fulfills the following equation
\begin{equation} \label{deq15}
 \widetilde{W}^\Psi(t) = \int_0^t A(t-\tau)\, \widetilde{W}^\Psi(\tau)\,d\tau
  + \int_0^t \Psi(\tau)\,dW(\tau)\,.
\end{equation}
\end{theor}

\proof{Let us notice that the process $\widetilde{W}^\Psi$ has integrable
trajectories.
Then, from Dirichlet's formula and stochastic Fubini's theorem
\begin{eqnarray*}
 \int_0^t A(t-\tau)\,\widetilde{W}^\Psi(\tau)\,d\tau \hspace{10ex} 
 &=& \int_0^t A(t-\tau)  \!
 \int_0^\tau \! \widetilde{S}(\tau-\sigma)\Psi(\sigma)dW(\sigma)\,d\tau = \\
 &=&  \int_0^t \!\left[ \int_\sigma^t \!
  A(t-\tau)\widetilde{S}(\tau-\sigma)\,d\tau\right]
  \Psi(\sigma)dW(\sigma) = \\  (z:= \tau-\sigma)\hspace{8ex}
 &=& \int_0^t \! \left[ \int_0^{t-\sigma} \!
  A(t-\sigma-z) \widetilde{S}(z)dz \right]
  \Psi(\sigma)dW(\sigma) = \\ 
 = \int_0^t (A\star \widetilde{S})(t-\sigma)\Psi(\sigma)dW(\sigma) 
 &=& \int_0^t [ \widetilde{S}(t-\sigma)-I] \Psi(\sigma)dW(\sigma) =\\
 = \!\int_0^t \!\! \widetilde{S}(t-\sigma)\Psi(\sigma)dW(\sigma) &-&
  \int_0^t \Psi(\sigma)dW(\sigma)\,.
\end{eqnarray*}
What gives the theorem.\hfill $\blacksquare$}

\section{Case $A(t)=a(t)A$}

In this section we consider the case when $A(t)=a(t)A$ and, in
contrary to the previous assumptions, $A$ is a closed {\tt unbounded}
operator in $H$ with a dense domain D(A). 
Then the equation (\ref{deq1}) has the form
\begin{equation} \label{deq15a}
 X(t) =  X_0 +\int_0^t a(t-\tau)A\,X(\tau)\,d\tau
  + \int_0^t \Psi(\tau)\,dW(\tau)\,, \quad t\in [0,T]\,.
\end{equation}
For details, precise setting and several results concerning (\ref{deq15a}) we
refer to (Karczewska, 2005). Because 
the results analogous to those from sections \ref{prob} and \ref{prop}
are formulated there, we do not repeat them here.

We denote 
$$ W^\Psi(t):= \int_0^t S(t-\tau)\Psi(\tau)dW(\tau), \quad t\in[0,T];
$$
now, the operators $S(t),~ t\geq 0$, become
resolvents (see (Pr\"uss 1993) for definition and properties).

In this case, the equality (\ref{deq15}) reads
\begin{equation} \label{deq16}
 W^\Psi(t) = \int_0^t a(t-\tau)A\, W^\Psi(\tau)\,d\tau
  + \int_0^t \Psi(\tau)\,dW(\tau)\,.
\end{equation}
Now, although the equation (\ref{deq15a}) is formally simpler than (\ref{deq1}),
the case is interesting because $A$ is unbounded. In consequence,
we are not able to obtain the equality (\ref{deq16}) 
directly like in Theorem \ref{th1}.

\begin{prop} \label{pr5}
If $\Psi\in\mathcal{N}^2(0,T;L_2^0)$ and $\Psi(\cdot,\cdot)(U_0)\subset D(A),$
 P-a.s., then the stochastic convolution 
$W^\Psi$ fulfills the equation
$$ \langle W^\Psi(t),\xi\rangle_H = \int_0^t \langle a(t-\tau)W^\Psi(\tau),
 A^*\xi\rangle_H + \int_0^t \langle \xi,\Psi(\tau)dW(\tau)\rangle_H,
\quad P-a.s.,
$$
for any $t\in [0,T]$ and $\xi\in D(A^*)$.
\end{prop}
{\bf Comment }{\em Proposition \ref{pr5} states that in case $A(t)=a(t)A$,
a mild solution to (\ref{deq15a}) is a weak solution to
(\ref{deq15a}).}\\

This is well-known (see, e.g.\ (Pr\"uss 1993)) that
the scalar equation corresponding to (\ref{deq15a}) is 
\begin{equation} \label{deq16a}
s(t) + \mu (a\star s)(t) = 1\,.
\end{equation}

\begin{defin} \label{df16}
 We say that function $a\in L^1(0,T)$ is {\tt completely positive}
 on  $[0,T]$ if for any $\mu\geq 0$, the solution to (\ref{deq16a})
 satisfies $s(t)\geq 0$ on $[0,T]$.
\end{defin}
\begin{prop} \label{pr5a} (Clem\'ent and Nohel 1979)
~Assume that $a\in L^1(0,T)$ and $a$ is nonnegative and nonincreasing 
on $[0,T]$. Then $a$ is comletely positive on $[0,T]$.
\end{prop}
{\bf Comment }{\em  If $a\in L^1(0,T)$ and is completely monotonic on 
$[0,T]$, i.e.\ $(-1)^ka^{(k)}(t)\geq 0$, $t\in (0,T),~k=0,1,\ldots$, then
$a$ is comletely positive on $[0,T]$.}

\begin{prop} \label{pr6} 
Assume that $A$ is m-accretive operator and function $a$ is completely
positive on $[0,T]$.
If $\Psi$ and $A\Psi$ belong to $\mathcal{N}^2(0,T;L_2^0)$ and in addition 
$\Psi(\cdot,\cdot)(U_0)\subset D(A),$ P-a.s., 
then the equality (\ref{deq16}) holds.
\end{prop}
\proof{The operator $A$ is m-accretive if an only if it generates 
the linear strongly continuous semigroup of contractions, see e.g.\ (Zheng
2004). 
So, we may use results due to (Clem\'ent and Nohel 1979).

Denote by $A_\lambda:=\frac{1}{\lambda}(I-J_\lambda)$, where 
$J_\lambda=(I+\lambda A)^{-1}$, $\lambda\geq 0$, the Yosida approximation of
the operator $A$. By $S_\lambda(t), t\geq 0$, we denote the resolvent
operators corresponding to the Volterra equation 
(\ref{deq15a}) with the operator 
$A_\lambda$ instead of the operator $A$. The paper (Clem\'ent and Nohel 1979) 
provides the convergence
$ \lim_{\lambda\rightarrow 0^+} S_\lambda(t)x = S(t)x \mbox{~~for~~}
 t\in [0,T] \mbox{~~and~~} x\in D(A).$

Let us recall, that  the formula (\ref{deq16})
holds for any bounded operator $A$. Then it holds for operators $A_\lambda$,
too: $$
W_\lambda^\Psi(t) = \int_0^t A_\lambda W_\lambda^\Psi(\tau)d\tau +
 \int_0^t\Psi(\tau)dW(\tau),$$ 
where $$ W_\lambda^\Psi(t) := \int_0^t S_\lambda(t-\tau)\Psi(\tau)dW(\tau).
$$

Because $A_\lambda x = AJ_\lambda x$ for $x\in D(A)$, then
$$ A_\lambda W_\lambda^\Psi(t) = 
 A_\lambda \int_0^t S_\lambda(t-\tau)\Psi(\tau)dW(\tau) =
 J_\lambda \int_0^t S_\lambda(t-\tau)A\Psi(\tau)dW(\tau).
$$

Now, from the properties of stochastic integral and by the Lebesgue dominated 
convergence theorem, we obtain $$
\lim_{\lambda\rightarrow 0^+} \sup_{t\in [0,T]} \mathbb{E}
  |W_\lambda^\Psi(t)-W^\Psi(t)|_H^2=0 $$
and $$
\lim_{\lambda\rightarrow 0^+} \sup_{t\in [0,T]} \mathbb{E}
  |A_\lambda W_\lambda^\Psi(t)-AW^\Psi(t)|_H^2=0 $$   
what gives the required result.  
\hfill $\blacksquare$}

\begin{defin}\label{df1a}
 Suppose $S(t),~t\geq 0$, is a resolvent for (\ref{deq15a}). 
 $S(t)$ is called
 {\tt expo\-nent\-ially bounded} if there are constants $M\geq 1$ and
 $w\in\mathbb{R}$ that $||S(t)||\leq M\,e^{wt}$ for all $t\geq 0$.
\end{defin}

In contrary to the case of semigroups, resolvents (if they exist) need
not to be exponentially bounded even if the kernel $a$ belongs to
$L^1(\mathbb{R}_+)$. Existence of such resolvents is given, e.g.\
by Theorem 1.3 in (Pr\"uss 1993). An important class of kernels 
providing exponentially bounded resolvent are 
$a(t)=t^{\alpha-1}/\Gamma(\alpha),~\alpha\in (0,2)$ or the class 
of completely monotonic functions.

\begin{prop} \label{pr7}
Assume that $A$ is the infinitesimal generator of a $C_0$--semi\-group
which is exponentially bounded and the kernel $a$ is a completely positive 
function. If $\Psi$ and $A\Psi$ belong to
$\mathcal{N}^2(0,T;L_2^0)$ and in addition 
$\Psi(\cdot,\cdot)(U_0)$ $\subset D(A),$ P-a.s., 
then the equality (\ref{deq16}) holds.
\end{prop}
The proof bases on the following auxiliary result.
\begin{lem} \label{lm2}
(Karczewska and Lizama 2005) Let $A$ be the infinitesimal generator of 
a $C_0$-semigroup 
$T(t),~t\geq 0$, satisfying $||T(t)||\leq Me^{\omega t}$ with
$M\geq 1,~\omega\in\mathbb{R}$
and suppose that 
$a\in L_\mathrm{loc}^1(\mathbb{R}_+)$ is a completely positive function.
Then $(A,a)$ admits an exponentially bounded resolvent $S(t)$.
Moreover, there exists bounded operators $A_n$, such that  $(A_n,a)$
admits resolvent family $S_n(t)$ satisfying 
$||S_n(t)||\leq Me^{w_0 t},~\omega_0\in\mathbb{R}$,
for all $t\geq 0,~n\in \mathbb{N}$, and
$ S_n(t) x \rightarrow S(t)x  $
for all $x\in X,~t\geq 0$. Moreover, the convergence is uniform in $t$ 
on every compact subset of $\mathbb{R}_+$.
\end{lem}

\end{document}